\documentclass[12pt] {article}
\usepackage{amsmath,amsthm}
\usepackage{amssymb,latexsym}
\title{Valuations on convex sets, non-commutative determinants, and pluripotential theory.}
\date{}
\author{ Semyon Alesker
\\  { \normalsize Department of Mathematics, Tel Aviv University, Ramat Aviv}
 \\  { \normalsize 69978 Tel Aviv,
Israel }
\\ {\normalsize e-mail: semyon@post.tau.ac.il}}
\def\RR{\mathbb{R}}
\def\CC{\mathbb{C}}

\def\ZZ{\mathbb{Z}}
\def\HH{\mathbb{H}}

\def\eps{\varepsilon}
\def\alp{\alpha}
\def\Ome{\Omega}
\def\lam{\lambda}

\def\d{det}
\def\qed { Q.E.D. }
\def\pr{\partial}

\def\dfq{\frac{\partial ^2 f}{\partial\bar q_i \partial q_j}}
\def\duq{\frac{\partial ^2 u}{\partial\bar q_i \partial q_j}}
\swapnumbers
\newtheorem{theorem}{Theorem}[subsection]
\newtheorem{corollary}[theorem]{Corollary}
\newtheorem{lemma}[theorem]{Lemma}
\newtheorem{proposition}[theorem]{Proposition}
\newtheorem{claim}[theorem]{Claim}
\theoremstyle{definition}
\newtheorem{example}[theorem]{Example}
\newtheorem{definition}[theorem]{Definition}
\newtheorem{remark}[theorem]{Remark}
\theoremstyle{proposition-definition}
\newtheorem{proposition-definition}[theorem]{Proposition-Definition}

\def\cf{{\cal F}}

  \def\cc{{\cal C}}
  \def\cf{{\cal F}}
 \def\ch{{\cal H}} 
 \def\ck{{\cal K}} 
  \def\co{{\cal O}}
\def\cp{{\cal P}}

\def\pt{\partial}
\def\ccw{{}\!^ {\textbf{C}} W}
\def\bcw{\bar{{}\!^{\textbf{C}} W}}
\def\chh{{}\!^ {\textbf{C}} \HH}
\def\rv{{}\!^ {\textbf{R}} V}
\begin{document}
\maketitle
\begin{abstract}
A new method of constructing translation invariant continuous
valuations on convex subsets of the quaternionic space $\HH^n$ is
presented. In particular new examples of $Sp(n)Sp(1)$-invariant
translation invariant continuous valuations are constructed. This
method is based on the theory of plurisubharmonic functions of
quaternionic variables developed by the author in two previous
papers \cite{alesker-bsm} and \cite{alesker-jga}.
\end{abstract}
\setcounter{section}{-1}
\section{Introduction.}\label{intro}
\setcounter{subsection}{1}

The goal of this paper is to present a new method of constructing
translation invariant continuous valuations on convex subsets of
the quaternionic space $\HH^n$. As an application of this method
we obtain new examples of $Sp(n)Sp(1)$-invariant translation
invariant continuous valuations. This method is based on the
theory of plurisubharmonic functions of quaternionic variables
developed by the author in two previous papers \cite{alesker-bsm}
and \cite{alesker-jga}. The main results of the paper are Theorem
\ref{quagen} and its immediate Corollary \ref{invariant}.


Let us remind basic notions of the theory of valuations on convex
sets referring for more details to the surveys by McMullen
\cite{mcmullen-survey} and McMullen and Schneider
\cite{mcmullen-schneider}. Let $V$ be a finite dimensional real
vector space. Let ${\cal K}(V)$ denote the class of all convex
compact subsets of $V$.

\begin{definition} a) A function $\phi :{\cal K}(V) \to \CC$ is called a
valuation if for any $K_1, \, K_2 \in {\cal K}(V)$ such that their
union is also convex one has
$$\phi(K_1 \cup K_2)= \phi(K_1) +\phi(K_2) -\phi(K_1 \cap K_2).$$

b) A valuation $\phi$ is called continuous if it is continuous
with respect the Hausdorff metric on ${\cal K}(V)$.
\end{definition}

Remind that the Hausdorff metric $d_H$ on ${\cal K}(V)$ depends on
a choice of a Euclidean metric on $V$ and it is defined as
follows: $d_H(A,B):=\inf\{ \eps >0|A\subset (B)_\eps \mbox{ and }
B\subset (A)_\eps\},$ where $(U)_\eps$ denotes the
$\eps$-neighborhood of a set $U$. Then ${\cal K} (V)$ becomes a
locally compact space, and the topology on $\ck(V)$ induced by the
Hausdorff metric does not depend on the choice of the Euclidean
metric on $V$.

The theory of valuations has numerous applications in convexity
and integral geometry (see e.g. \cite{alesker-jdg},
\cite{hadwiger-book}, \cite{klain-rota}, \cite{schneider-book}).
Classification results on valuations invariant under specific
groups are of particular importance for applications. The famous
result of H. Hadwiger \cite{hadwiger-book} describes explicitly
isometry invariant continuous valuations on a Euclidean space
$\RR^n$. For discussion of applications of Hadwiger's result to
integral geometry in Euclidean spaces we refer to the books
Hadwiger \cite{hadwiger-book}, Klain and Rota \cite{klain-rota},
and the survey article by Hug and Schneider \cite{hug-schneider}.

It was shown by the author in \cite{alesker-adv} that if $G$ is a
compact subgroup of the orthogonal group $O(n)$ acting
transitively on the unit sphere then the space of $G$-invariant
translation invariant continuous valuations (let us denote this
space by $Val^G$) is finite dimensional. Hence in this case one
may hope to obtain an explicit description of the space $Val^G$.
It was shown in \cite{alesker-mult} and \cite{alesker-lefschetz}
(see also \cite{alesker-jdg}) that for a group $G$ satisfying the
above assumptions the space $Val^G$ has various remarkable
properties: it has a canonical structure of a commutative
associative graded algebra satisfying the Poincar\'e duality, and
if $-Id\in G$ it satisfies the hard Lefschetz type theorem. Remind
also that there is an explicit classification of compact connected
groups acting transitively and effectively on the sphere
(\cite{borel1}, \cite{borel2}, \cite{montgomery-samelson}). Namely
there are 6 infinite series $SO(n), U(n), SU(n), Sp(n),
Sp(n)Sp(1), Sp(n)U(1)$, and 3 exceptions $G_2, Spin(7), Spin(9)$.

The case of $Val^G$ for the groups $G=O(n)$ and $G=SO(n)$ was
completely classified by Hadwiger \cite{hadwiger-book}. In
\cite{alesker-jdg} the author has obtained an explicit
classification of translation invariant $U(n)$-invariant
continuous valuations on an Hermitian space $\CC^n$; in that
article also some applications to integral geometry in Hermitian
spaces were obtained. The space of $SU(2)$-invariant valuations on
$\CC^2$ was described in \cite{alesker-su2}.

In this article we study the case of translation invariant
$Sp(n)Sp(1)$-invariant valuations on the quaternionic space
$\HH^n$. Though we do not obtain a complete classification, we
construct new non-trivial examples of such valuations using a new
method. This method  is based on the theory of plurisubharmonic
functions of quaternionic variables developed by the author in
\cite{alesker-bsm} and \cite{alesker-jga}. Motivated by
applications to valuations, we obtain in this paper further
results in this theory. In order to explain our main results let
us first discuss their complex analogs which are more classical
and well known.

Let $\CC^n$ be a Hermitian space with the Hermitian product
$(\cdot, \cdot)$. For a convex compact subset $K\in \ck(\CC^n)$
let $h_K$ denote its supporting functional (remind that
$h_K:\CC^n\to \RR$, and $h_K(x):=\sup_{y\in K}(x,y)$). Let us fix
an integer $1\leq l\leq n$. Let us fix a continuous compactly
supported $(n-l,n-l)$-form $\psi$ on $\CC^n$. Define
\begin{eqnarray}\label{intr-ex}
\phi(K):=\int_{\CC^n}(dd^ch_K)^l\wedge \psi.
\end{eqnarray}
Theorem \ref{kazgen} of this paper claims that $\phi$ is a
translation invariant continuous valuation. Note that since the
function $h_K$ is not necessarily smooth (but it is convex, and
hence continuous plurisubharmonic) the expression $(dd^c h_K)^l$
should be understood in the sense of currents using the
Chern-Levine-Nirenberg theorem \cite{chern-levine-nirenberg}. The
continuity property of $\phi$ is not obvious and also follows from
the same result \cite{chern-levine-nirenberg}. The property of
valuation is not evident; it is a consequence of a more general
result about plurisubharmonic functions due to Z. Blocki
\cite{blocki} (this was kindly explained to us by N. Levenberg).

The main result of this article is a construction of a
quaternionic analogue on $\HH^n$ of valuations of the form
(\ref{intr-ex}) (Theorem \ref{quagen}). This construction is based
on the theory of plurisubharmonic (psh) functions of quaternionic
variables. The notion of psh function of quaternionic variables
was introduced by the author in \cite{alesker-bsm} and
independently by G. Henkin \cite{henkin}. This class of functions
was studied further by the author in \cite{alesker-jga} where the
quaternionic Monge-Amp\`ere equations were introduced and
investigated. A quaternionic version of the Chern-Levine-Nirenberg
theorem was proved in \cite{alesker-bsm}. In this paper we prove a
refined version of that result (Theorem \ref{cont}). (The
refinement is approximately as follows: our previous result from
\cite{alesker-bsm} corresponds in the complex case to establishing
of some properties of the current $(dd^ch)^n$ on $\CC^n$ where $h$
is a continuous complex psh function, and Theorem \ref{cont} of
this paper corresponds to the current $(dd^c h)^l$ with $0\leq
l\leq n$.) This required to introduce a quaternionic analogue of
the notion of positive current (Subsection \ref{poscur}). Then we
prove a quaternionic analogue of Blocki's formula (Theorem
\ref{blocki}).

Note that most of the results of this article and
\cite{alesker-bsm}, \cite{alesker-jga} as well make use of the
notion of non-commutative determinant, particularly the Moore
determinant of quaternionic hyperhermitian matrices which is
reviewed in Subsection \ref{det}. For the purposes of this article
it was necessary to present a new, coordinate free, construction
of the Moore determinant. This is done in Subsection \ref{forms}.
Moreover we have constructed a quaternionic analogue of the
algebra of exterior forms of type $(p,p)$ on a complex space
(Definition \ref{forms-12}).
Note also that it was shown in \cite{gelfand-retakh-wilson} that
the Moore determinant can be expressed via Gelfand-Retakh
quasideterminants first introduced in \cite{gelfand-retakh} and
which generalize most of the known notions of non-commutative
determinants (see also Remark \ref{remark-deter} in Subsection 1.2
of this paper). For the details we refer to the recent survey
\cite{GGRW}, and for more details on the quaternionic case we
refer to \cite{gelfand-retakh-wilson}.

Let us describe briefly the content of the paper. Section
\ref{back} does not contain new results. In Subsection \ref{ca} we
review some relevant facts from the complex analysis. In
Subsection \ref{det} we remind the necessary definitions and facts
about hyperhermitian matrices and the Moore determinant. In
Subsection \ref{qpt} we summarize the relevant definitions and
results about psh functions of quaternionic variables following
\cite{alesker-bsm}.

Section \ref{la} contains some new constructions from quaternionic
linear algebra. Subsection \ref{forms} is purely algebraic. There
we describe a coordinate free construction of the Moore
determinant. Next for a right quaternionic $\HH$-module $V$ of
finite dimension we construct a graded algebra $\Omega^\bullet
(V)$ analogous to the graded algebra of exterior forms of type
$(p,p)$ on a complex space. In Subsection \ref{posforms} we
introduce the notions of weakly and strongly positive elements in
$\Omega^\bullet (V)$.

In Section \ref{psh} we study further psh functions of
quaternionic variables. In Subsection \ref{poscur} we prove a
quaternionic analogue of the Chern-Levine-Nirenberg theorem
refining our previous result from \cite{alesker-bsm}. This theorem
concerns positive currents with values in $\Omega^\bullet (V)$. In
Subsection \ref{qblocki} we obtain a quaternionic analogue of the
Blocki formula.

Section \ref{valuations} contains applications of the above
results to the theory of valuations. In Subsection \ref{kazar} we
remind Kazarnovskii's pseudovolume and its generalizations
following \cite{kazarnovskii-81} and \cite{kazarnovskii-84}. In
Subsection \ref{vq} we obtain quaternionic analogues of these
valuations, i.e. the main results of this paper (Theorem 4.2.1 and
Corollary 4.2.2).

{\bf Acknowledgements.} We express our gratitude to N. Levenberg
who has informed us about Blocki's paper \cite{blocki} and has
provided an argument to deduce from it Theorem \ref{kazgen}. We
thank also I.M. Gelfand, G. Henkin, A. Rashkovskii, V. Retakh, and
M. Sodin for useful conversations.

\section{Background.}\label{back} In this section we remind some
definitions and facts from complex analysis  which are very
classical. Then we remind some facts from the theory of
non-commutative determinants, the exposition will follows
essentially \cite{alesker-bsm}. Then we review  the theory of
plurisubharmonic functions of quaternionic variables developed by
the author in \cite{alesker-bsm} and which is based on the theory
of non-commutative determinants. This section does not contain new
results.
\subsection{Some complex analysis.}\label{ca} Let us remind the definition
of a plurisubharmonic function of complex variables (see e.g.
Lelong's book \cite{lelong} for more details). Let $\Ome$ be an
open subset in $\CC^n$.
\begin{definition}\label{ca-1}
A real valued function $u: \Omega \to \RR$ is called
plurisubharmonic if it is upper semi-continuous and its
restriction to any {\itshape complex} line is subharmonic.
\end{definition}
 Recall that upper semi-continuity means that
 $u(x_0)\geq \underset{x\to x_0}{\limsup u(x)}$ for any $x_0\in \Omega$.
 We will denote by $P(\Omega)$ the class of plurisubharmonic (psh)
 functions in the domain $\Omega$,
and the class of continuous functions in $\Omega$ will be denoted
by $C(\Ome)$.
\begin{remark}\label{ca-2}
If in the above definition one replaces the word "complex" by the
word "real" everywhere then one obtains a definition equivalent to
the usual definition of {\itshape convex} function.
\end{remark}
The following result is due to Chern, Levine, and Nirenberg
\cite{chern-levine-nirenberg} (its real analogue is due to A.D.
Aleksandrov \cite{aleksandrov-58}).
\begin{theorem}[\cite{chern-levine-nirenberg}]\label{ca-3}
Let $\Ome \subset \CC^n$ be an open subset. For any function $u\in
C(\Omega) \cap P(\Omega)$ one can define a non-negative measure
denoted by $\det (\frac{\pt ^2 u}{\pr\bar z_i\pt  z_j})$ which is
uniquely characterized by the following two properties:
\newline
(a) if $u\in C^2(\Omega)$ then it has the obvious meaning;
\newline
(b) if $u_N\to u$ uniformly on compact subsets in $\Ome$, and
$u_N,\, u\in C(\Omega)\cap P(\Omega)$, then
$$\det (\frac{\partial ^2 u_N}{\partial\bar z_i \partial z_j})\overset
{w}{\to} \det(\frac{\pt ^2 u}{\pr\bar z_i\pt z_j}),$$ where the
convergence of measures in understood in the sense of weak
convergence of measures.
\end{theorem}
\begin{remark}\label{ca-4}
In fact it is easy to see from the definition that the uniform
limit of functions from $P(\Ome)\cap C(\Ome)$ also belongs to this
class.
\end{remark}
\subsection{Non-commutative determinants.}\label{det} For our purposes we
will need the notions of Moore determinant. The exposition follows
essentially \cite{alesker-bsm}. For more developed theory of
non-commutative determinants, so called quasideterminants of
Gelfand-Retakh, we refer to the survey \cite{GGRW}. In
\cite{gelfand-retakh-wilson} there was established the connection
of quasideterminants to the Moore determinant. A good survey of
quaternionic determinants is \cite{aslaksen}.

\begin{definition}\label{det-5}
  A {\itshape hyperhermitian semilinear form}
on $V$ is a map $ a:V \times V \to \HH$ satisfying the following
properties:

(a) $a$ is additive with respect to each argument;

(b) $a(x,y \cdot q)= a(x, y) \cdot q$ for any $x,y \in V$ and any
$q\in \HH$;

(c) $a(x,y)= \overline{a(y,x)}$.
\end{definition}

\begin{example}\label{det-6} Let $V= \HH ^n$ be the standard coordinate space
considered as right vector space over $\HH$. Fix a {\itshape
hyperhermitian} $n \times n$-matrix $(a_{ij})_{i,j=1}^{n}$, i.e.
$a_{ij} =\bar a_{ji}$, where $\bar x$ denotes the usual
quaternionic conjugation of $x$. For $x=(x_1, \dots, x_n), \,
y=(y_1, \dots, y_n)$ define
$$A(x,y) = \sum _{i,j} \bar x_i a_{ij} y_j$$
(notice the order of the terms!). Then $A$ defines hyperhermitian
semilinear form on $V$.

The set of all hyperhermitian $n\times n$-matrices will be denoted
by $\ch_n$. Then $\ch_n$ a vector space over $\RR$.
\end{example}
In general one has the following standard claims.

\begin{claim}\label{det-7} Fix a basis in a finite dimensional right quaternionic
vector space $V$. Then there is a natural bijection between the
space of hyperhermitian semilinear forms on $V$ and the space
$\ch_n$ of $n \times n$-hyperhermitian matrices.
\end{claim}
This bijection is in fact described in previous Example 1.2.2.

\begin{claim}\label{det-8} Let $A$ be the matrix
of the given hyperhermitian form in the given basis. Let $C$ be
transition matrix from this basis to another one. Then the matrix
$A'$ of the given form in the new basis is equal $$A' =C^* AC .$$
\end{claim}
\begin{remark}\label{det-9}
Note that for any hyperhermitian matrix $A$ and for any matrix $C$
the matrix $C^* AC$ is also hyperhermitian. In particular the
matrix $C^* C$ is always hyperhermitian.
\end{remark}

\begin{definition}\label{det-10} A hyperhermitian semilinear form $a$
is called {\itshape positive definite} if $a(x,x)>0$ for any
non-zero vector $x$. Similarly $a$ is called {\itshape
non-negative definite} if $a(x,x)\geq 0$ for any vector $x$.
 \end{definition}

Let us fix on our quaternionic right vector space $V$ a positive
definite hyperhermitian form $( \cdot , \cdot )$. The space with
fixed such a form will be called {\itshape hyperhermitian} space.

For any quaternionic linear operator $\phi: V\to V$ in
hyperhermitian space one can define the adjoint operator $\phi ^*
:V \to V$ in the usual way, i.e. $(\phi x,y)= (x, \phi ^* y)$ for
any $x,y \in V$. Then if one fixes an orthonormal basis in the
space $V$ then the operator $\phi$ is selfadjoint if and only if
its matrix in this basis is hyperhermitian.

\begin{claim}\label{det-11}
For any  selfadjoint operator in a hyperhermitian space there
exists an orthonormal basis such that its matrix in this basis is
diagonal and real.
\end{claim}
 Now we are going to define the Moore
determinant of  hyperhermitian matrices. The definition below is
different from the original one \cite{moore} but equivalent to it.

Any quaternionic matrix $A\in M_n(\HH)$ can be considered as a
matrix of an $\HH$-linear endomorphism of $\HH^n$. Identifying
$\HH^n$ with $\RR^{4n}$ in the standard way we get an $\RR$-linear
endomorphism of $\RR^{4n}$. Its matrix in the standard basis will
be denoted by ${}^{\mathbb{R}} A$, and it is called the
realization of $A$. Thus ${}^{\mathbb{R}} A\in M_{4n}(\RR)$.

 Let us consider the entries of $A$ as formal variables (each
quaternionic entry corresponds to four commuting real variables).
Then $det ({}^{\mathbb{R}} A)$  is a homogeneous polynomial of
degree $4n$ in $n(2n-1)$ real variables. Let us denote by $Id$ the
identity matrix.
 One has the following result.
\begin{theorem}\label{det-12}
There exists a polynomial $P$ defined on the space $\ch_n$ of all
hyperhermitian $n \times n$-matrices such that for any
hyperhermitian $n \times n$-matrix $A$ one has
$det({}^{\mathbb{R}} A)= P^4(A)$ and $P(Id)=1$. $P$ is defined
uniquely by these two properties. Furthermore $P$ is homogeneous
of degree $n$ and has integer coefficients.
\end{theorem}
Thus for any hyperhermitian matrix $A$ the value $P(A)$ is a real
number, and it is called the {\itshape Moore determinant} of the
matrix $A$. The explicit formula for the Moore determinant  was
given by Moore \cite{moore} (see also \cite{aslaksen}). From now
on the Moore determinant of a matrix $A$ will be denoted by $det
A$. This notation should not cause any confusion with the usual
determinant of real or complex matrices due to part (i) of the
next theorem.
\begin{theorem}\label{det-13}

(i) The Moore determinant of any complex hermitian matrix
considered as quaternionic hyperhermitian matrix is equal to its
usual determinant.

(ii) For any hyperhermitian $n\times n$-matrix $A$ and any matrix
$C\in M_n(\HH)$
$$det (C^*AC)= detA \cdot det(C^*C).$$
\end{theorem}
\begin{example}\label{det-14}

(a) Let $A =diag(\lam_1, \dots, \lam _n)$ be a diagonal matrix
with real $\lam _i$'s. Then $A$ is hyperhermitian and the Moore
determinant $detA= \prod _{i=1}^n \lam_i$.

(b)  A general hyperhermitian $2 \times 2$-matrix $A$ has the form
 $$ A=  \left[ \begin {array}{cc}
                     a&q\\
                \bar q&b\\
                \end{array} \right] ,$$
where $a,b \in \RR, \, q \in \HH$. Then $det A =ab - q \bar q$.
\end{example}

\begin{definition}\label{det-14.5}
A hyperhermitian $n\times n$-matrix $A=(a_{ij})$ is called
{\itshape positive} (resp. {\itshape non-negative}) {\itshape
definite} if for any non-zero vector $\xi= \left[\begin{array}{c}
             \xi_1\\
             \vdots\\
             \xi_n
             \end{array}\right]$ one has
$\xi^*A\xi=\sum_{ij}\bar\xi_ia_{ij}\xi_j
>0$ (resp. $\geq 0$).
\end{definition}

\begin{claim}\label{det-15}
Let $A$ be  a non-negative (resp. positive) definite
hyperhermitian matrix. Then $det A \geq 0 \, (\mbox{ resp. } det A
>0)$.
\end{claim}

Let us remind now the definition of the mixed discriminant of
hyperhermitian matrices in analogy with the case of real symmetric
matrices \cite{aleksandrov-38}.
\begin{definition}\label{det-16}
Let $A_1, \dots ,A_n$ be hyperhermitian $n \times n$- matrices.
Consider the homogeneous polynomial in real variables $\lam _1
,\dots , \lam _n$ of degree $n$ equal to $det(\lam_1 A_1 + \dots +
\lam_n A_n)$. The coefficient of the monomial $\lam_1 \cdot \dots
\cdot \lam_n$ divided by $n!$ is called the {\itshape mixed
discriminant} of the matrices  $A_1, \dots ,A_n$, and it is
denoted by $\d(A_1, \dots ,A_n)$.
\end{definition}
Note that the mixed discriminant is symmetric with respect to all
variables, and linear with respect to each of them, i.e.
$$\d (\lam A_1' +\mu A_1'', A_2, \dots, A_n )=
\lam \cdot \d( A_1', A_2, \dots, A_n ) + \mu \cdot \d(A_1'', A_2,
\dots, A_n )$$ for any {\itshape real} $\lam , \, \mu$. Note also
that $\d(A, \dots, A)=det A$.

\begin{theorem}\label{det-17}
 The mixed discriminant of positive (resp. non-negative)
definite matrices is positive (resp. non-negative).
\end{theorem}


We will need also the following proposition which is essentially
well known but we do not have an exact reference.
\begin{proposition}\label{det-18}
Let $A$ be a quaternionic $n\times n$-matrix. Then
$$\det({}^{\mathbb{R}} A)=(\det(A^*A))^2=(\det (AA^*))^2$$
where the first determinant denotes the usual determinant of the
real $(4n\times 4n)$-matrix ${}^{\mathbb{R}} A$, and the second
and the third determinants denote the Moore determinant.
\end{proposition}
{\bf Proof.} We may assume that $A$ is invertible, i.e. $A\in
GL_n(\HH)$ (otherwise both sides are equal to 0). The maximal
compact subgroup of $GL_n(\HH)$ is $Sp(n)$. There exists a
decomposition $A=U\cdot D\cdot V$ where $U,\, V\in Sp(n)$ and $D$
is a real diagonal matrix, $D= diag(\lam_1,\dots, \lam_n),\,
\lam_i\in \RR$. Since $Sp(n)\subset SO(4n)$ we have $$\det
({}^{\mathbb{R}} A)= \det({}^{\mathbb{R}}
D)=(\prod_{i=1}^n\lam_i)^4.$$ On the other hand
$$\det (AA^*)=\det(UD^2 U^*)=(\prod_{i=1}^n\lam_i)^2$$
where the last equality follows from Theorem \ref{det-13} (ii) and
Example \ref{det-14}. Hence $\det ({}^{\mathbb{R}}
A)=(\det(AA^*))^2$. Similarly one shows that $\det({}^{\mathbb{R}}
A)=(\det(A^*A))^2$. \qed

Let us introduce more notation. Let $A$ be any hyperhermitian $n
\times n$- matrix. For any subset $I \subset \{1, \dots, n\}$ the
minor $M_I(A)$ of $A$ which is obtained by deleting the rows and
columns with indexes from the set $I$, is clearly hyperhermitian.
(For instance $M_\emptyset (A)=A$.) For $I=\{1, \dots , n \}$ set
$\det M_{\{1, \dots,n \}} :=1$. We will need a lemma.
\begin{lemma}\label{det-19}
Let $T:=\left[ \begin{array}{ccc}
                     t_1&        &    0 \\
                        & \ddots & \\
                      0 &        &t_n \\
            \end{array} \right]$ be a real diagonal $n\times n$-matrix.
Let $A_1,\dots, A_{n-1}$ be hyperhermitian $n\times n$-matrices.
Then $$\det(T,A_1,\dots,A_{n-1})= \frac{1}{n}\sum_{i=1}^n t_i
\det(M_{\{i\}}(A_1),\dots,M_{\{i\}}(A_{n-1})).$$
\end{lemma}
{\bf Proof.} In \cite{alesker-bsm}, Proposition 1.1.11, it was
shown that for a matrix $T$ as above and for any hyperhermitian
$n\times n$-matrix $A$ one has
\begin{eqnarray}\label{det-form}
\det (A+T) = \sum _{I\subset \{1, \dots ,n\} } (\prod _{i\in I}
t_i) \cdot \det M_I(A).\end{eqnarray}
 Lemma \ref{det-19} follows from this formula
and the definition of the mixed discriminant. \qed

The following lemma also will be used later.
\begin{lemma}\label{det-20}
Let $A_1,\dots,A_n\geq 0$ be hyperhermitian $n\times n$-matrices.
Then
$$\det(A_1,\dots,A_n)\leq \det (\sum_{i=1}^nA_i).$$
\end{lemma}
{\bf Proof.} First observe that Theorem \ref{det-17} implies that
if $B_1,\dots,B_n,C_1,\dots,C_n$ are hyperhermitian $(n\times
n)$-matrices and $$B_i\leq C_i \mbox{ for } i=1,\dots, n$$ then
$$\det(B_1,\dots,B_n)\leq \det(C_1,\dots,C_n).$$ Now let us take
$B_i:=A_i$ and $C_1=\dots=C_n:=\sum_{i=1}^nA_i$. Lemma
\ref{det-20} is proved. \qed

\begin{remark}\label{remark-deter}
The Moore determinant of hyperhermitian quaternionic matrices
behaves exactly as the usual determinant of complex hermitian or
real symmetric matrices from all points of view. There is also a
notion of the Dieudonn\'e determinant \cite{dieudonne} of an
arbitrary quaternionic matrix which behaves exactly like the
absolute value of the usual determinant of complex or real
matrices. It was applied to the theory of quaternionic psh
functions in author's paper \cite{alesker-bsm}. Implicitly the
Dieudonn\'e determinant is also used in this article since the
proof in \cite{alesker-bsm} of Theorem \ref{qpt-2}  was based on
it. The Dieudonn\'e determinant can also be expressed via
Gelfand-Retakh quasideterminants (see \cite{GGRW}).
\end{remark}

\subsection{Quaternionic pluripotential theory.}\label{qpt} The notion of plurisubharmonic
(psh) function of quaternionic variables was introduced by the
author in \cite{alesker-bsm} and independently by G. Henkin
\cite{henkin}. The exposition here follows \cite{alesker-bsm}
where the definition of psh function of quaternionic variables is
presented in the form suggested by G. Henkin, and parallel to the
complex case. Let $\Omega $ be a domain in $\HH ^n$.
\begin{definition}\label{qpt-1}
A real valued function $u: \Omega \to \RR$ is called quaternionic
plurisubharmonic if it is upper semi-continuous and its
restriction to any right {\itshape quaternionic} line is
subharmonic.
\end{definition}

The class of psh functions in $\Ome$ will be denoted by $P(\Ome)$,
and the class of continuous functions will be denoted by
$C(\Ome)$.

Let $q$ be a quaternionic coordinate,
$$q=t+ix+jy +kz ,$$
where $t,x,y,z$ are real numbers. Consider the following operators
defined on the class of smooth $\HH$-valued functions of the
variable $q\in \HH$:
\def\db{\frac{\partial}{\partial \bar q}}
\def\dq{\frac{\partial}{\partial  q}}
$$\db f:=\frac{\partial f}{\partial  t}  +
i \frac{\partial f}{\partial  x} + j \frac{\partial f}{\partial y}
+ k \frac{\partial f}{\partial  z}  , \mbox { and }$$
$$\dq f:=\overline{ \db \bar f}=
\frac{\partial f}{\partial  t}  -
 \frac{\partial f}{\partial x}  i-
 \frac{\partial f}{\partial  y} j-
\frac{\partial f}{\partial  z}  k.$$ Note that $\db$ is called
sometimes Cauchy-Riemann-Moisil-Fueter operator, or sometimes
Dirac-Weyl operator, or just Dirac operator. It is easy to see
that $\db$ and $\dq$ commute, and if $f$ is a {\itshape real
valued} function then
$$\db \dq f= \Delta f = (\frac{\partial ^2 }{\partial t ^2}+
\frac{\partial ^2 }{\partial x ^2}+ \frac{\partial ^2 }{\partial y
^2}+ \frac{\partial ^2 }{\partial z ^2})f.$$ For any real valued
$C^2$- smooth function $f$ the matrix $(\dfq)_{i,j=1}^{n}$ is
obviously hyperhermitian.  Note also that the operators
$\frac{\partial}{\partial q_i}$ and $\frac{\partial}{\partial \bar
q_j}$ commute. One can easily check the following identities.

\begin{proposition}\label{qpt-1.5}
 Let $f:\HH ^n \to \HH$ be a smooth function. Then for any
$\HH$-linear transformation $A$ of $\HH ^n$ (as right $\HH
$-vector space) one has the identities
$$ \left( \frac {\partial ^2 f(Aq)}{\partial \bar q_i \partial q_j} \right)
=A^* \left(\frac {\partial ^2 f}{\partial \bar q_i \partial
q_j}(Aq) \right)A .$$
\end{proposition}

We have the following quaternionic analogue of Theorem \ref{ca-3}.
\begin{theorem}[\cite{alesker-bsm}]\label{qpt-2}
For any function $u\in C(\Omega) \cap P(\Omega)$ one can define a
non-negative measure  denoted by $\det (\duq)$ which is uniquely
characterized by the following two properties:
\newline
(a) if $u\in C^2(\Omega)$ then it has the obvious meaning;
\newline
(b) if $u_N\to u$ uniformly on compact subsets in $\Omega$, and
$u_N,\, u\in C(\Omega)\cap P(\Omega)$, then
$$\det (\frac{\partial ^2 u_N}{\partial\bar q_i \partial q_j})\overset
{w}{\to} \det(\duq),$$ where the convergence of measures is
understood in the sense of weak convergence of measures.
\end{theorem}
\begin{remark}\label{qpt-3}
In fact it is easy to see from the definition that the uniform
limit of functions from $P(\Ome)\cap C(\Ome)$ also belongs to this
class.
\end{remark}

\section{More quaternionic linear algebra.}\label{la}

\subsection{The space of forms.}\label{forms}
 The goal of this
subsection is to construct a quaternionic analogue of the spaces
of complex forms of type $(k,k)$. We will also present another
(coordinate free) construction of the Moore determinant.

\def\baw{\bar W}
Let $W$ be a right finite dimensional $\HH$-module,
$n=\dim_{\HH}W$. Let $\baw$ denote the quaternionic conjugate
space of $W$. Recall that $\baw$ is a left $\HH$-module, it
coincides with $W$ as a group, and the multiplication by scalars
from $\HH$ is given by $q\cdot w := w \cdot \bar q$. Consider the
tensor product $W\otimes _{\HH} \baw$. On this space one has an
involution $\sigma:W\otimes _{\HH} \baw\to W\otimes _{\HH} \baw $
defined by $\sigma(x\otimes y)=y\otimes x$.
\begin{theorem}\label{forms-1}
Let $W$ be a right $\HH$-module, $\dim_{\HH}W=n$. Then the space
$Sym^n((W\otimes_{\HH}\bar W)^\sigma)$ has unique
$Aut_{\HH}W$-invariant subspace of (real) codimension 1.
\end{theorem}
We postpone the proof of this theorem.
\begin{definition}\label{forms-2}
The one dimensional quotient space of $Sym^n(W\otimes_{\HH}\bar
W)$ by the subspace from Theorem \ref{forms-1} will be denoted by
$M(W)$. The quotient map $M: Sym^n((W\otimes_{\HH}\bar
W)^\sigma)\to M(W)$ will be called the {\itshape Moore map}.
\end{definition}
Let us explain why we call the map $M$ the Moore map. In fact it
coincides with the Moore determinant in the following sense.
Assume we are given on $W$ a hyperhermitian positive definite
form. Then it identifies the dual space
$W^*:=Hom_{\RR}(W,\RR)=Hom_{\HH}(W,\HH)$ with $\bar W$. Then
$W\otimes_{\HH}\bar W\simeq W\otimes _{\HH}W^*= End_{\HH}(W,W)$.
Under this identification $(W\otimes_{\HH}\bar W)^\sigma$
corresponds to selfadjoint endomorphisms of $W$. For any
selfadjoint operator $A$, $M(A,\dots, A)\in M(W)$. The space
$M(W)$ can be identified with $\RR$ as $Sp(n)$-module (since the
last group is compact and connected, and hence does not have
non-trivial real valued characters) when $1\in \RR$ corresponds to
$M(Id,\dots, Id)$. Under this identifications $M(A,\dots,A)$ is
equal to the Moore determinant of $A$.

Before we prove Theorem \ref{forms-1} we need some preparations.
Let us denote by $\chh$ the algebra of complex quaternions:
$\chh:=\HH\otimes _{\RR}\CC$. Recall that there exists an
isomorphism of algebras $\chh\simeq M_2(\CC)$. Denote also
$\ccw:=W\otimes_{\RR}\CC$, and $\bcw :=\bar W\otimes_{\RR}\CC$.
Let us fix a simple non-trivial right $\chh$-module $T$ (which is
unique up to non-canonical isomorphism since $\chh$ is a central
simple algebra). Then $\dim_{\CC}T= 2$.
\begin{remark}\label{forms-3.5}
Any finitely generated $\chh$-module is isomorphic to a direct sum
of finitely many copies of $T$. This is a general fact for modules
over central simple algebras (see e.g. \cite{weil}, Ch. IX \S 1,
Prop. 1).
\end{remark}

\begin{lemma}\label{forms-4}
There exists an isomorphism of $\CC$-vector spaces
$$(W\otimes_{\HH}\baw)^{\sigma}\otimes_{\RR}\CC=\wedge^2(Hom_{\chh}(T, \ccw))$$
which commutes with the action of $Aut_{\HH}W$.
\end{lemma}

\begin{remark}\label{forms-5}
(1) This isomorphism commutes with the natural action of
$(Aut_{\HH}W)\otimes_{\RR} \CC\simeq GL_{2n}(\CC)$.

(2) $\dim_{\CC}Hom_{\chh}(T, \ccw)=2n$. Indeed by Remark
\ref{forms-3.5} the $\chh$-module $\ccw$ is isomorphic to a direct
sum of $2n$ copies of the module $T$. Also one has isomorphism of
algebras $End_{\chh}(T,T)=\CC$.
\end{remark}
{\bf Proof} of Lemma \ref{forms-4}. Note that the vector space in
the left hand side is canonically isomorphic to
$(\ccw\otimes_{\chh}\bcw)^\sigma$. Denote $Z:=Hom_{\chh}(T,
\ccw)$. Then the evaluation map $T\otimes _{\CC}Z\to \ccw$ is an
isomorphism. Then $\bcw$ can be identified with $\bar T\otimes
_{\CC} Z$. Then $\ccw \otimes _{\chh} \bcw =Z\otimes_{\CC}(
T\otimes _{\chh}\bar T)\otimes _{\CC} Z $. But $T\otimes
_{\chh}\bar T$ is one dimensional. Hence $\ccw \otimes _{\chh}
\bcw$ is isomorphic to $Z\otimes_{\CC} Z$ and this identification
commutes with the action of $(Aut_{\HH}W)\otimes _{\RR}\CC=
Aut_{\CC}Z$. But $Z\otimes_{\CC}Z$ decomposes uniquely under the
action of $GL(Z)$ into two irreducible subspaces: $$Z\otimes Z=
\wedge^2 Z\oplus Sym^2 Z.$$ Hence on $Z\otimes Z$ there are only
two non-trivial involutions commuting with the action of $Aut_\CC
Z$. It is easy to see that $\sigma$ corresponds to that acting
trivially on $\wedge ^2 Z$ and by multiplication by $-1$ on $Sym^2
Z$. \qed

Now Theorem \ref{forms-1} follows from the next proposition.
\begin{proposition}\label{forms-6}
Let $Z$ be a complex $2n$-dimensional space. The the space
$Sym^n(\wedge^2 Z)$ contains a unique $GL(Z)$-invariant subspace
of codimension 1.
\end{proposition}
{\bf Proof.} First let us show the existence. Note that there
exists a canonical map $Sym^n(\wedge^2 Z)\to \wedge^{2n}Z$ given
by $x_1\otimes \dots \otimes x_n\mapsto x_1\wedge\dots \wedge
x_n$. Since $\wedge^{2n}Z$ is one dimensional the existence
follows.

The representation of $GL(Z)$ in $Sym^n(\wedge^2 Z)$ is algebraic
and hence completely reducible. Moreover the center of $GL(Z)$
acts by a multiplication by a character. Hence all one dimensional
components in $Sym^n(\wedge^2 Z)$ must be isomorphic. Thus the
uniqueness follows from the following result due to Howe
\cite{howe} which also will be used later.
\begin{proposition}[\cite{howe}, p. 563, Proposition
2]\label{forms-7} Let $Z$ be an even dimensional complex vector
space. Let $k\geq 0$ be an integer. The natural representation of
$GL(Z)$ in $Sym^k(\wedge^2 Z)$ is multiplicity free.
\end{proposition}
Thus Theorem \ref{forms-1} is proved. \qed

Now we will introduce a quaternionic analogue of the space of
complex vectors of type $(k,k)$. Let $W$ be a right $\HH$-module,
$\dim_{\HH}W=n$. The space $Sym^\bullet (W\otimes_{\HH}\bar
W)^\sigma$ is a commutative associative  algebra (it is full
symmetric algebra of a real vector space). For $0\leq k\leq n$ one
has a map
$$Sym^k(W\otimes_{\HH}\bar
W)^\sigma \times Sym^{n-k}(W\otimes_{\HH}\bar W)^\sigma\to M(W)$$
which is a composition of the product with the Moore map $M$. Let
us denote by $L^k(W)$ the left kernel of this map, namely
$$L^k(W)=\{x\in Sym^k(W\otimes_{\HH}\bar W)^\sigma |\, M(x\cdot
y)=0 \forall y\in Sym^{n-k}(W\otimes_{\HH}\bar W)^\sigma \}.$$ Let
us define $R^k(W):=Sym^k(W\otimes_{\HH}\bar W)^\sigma /L^k(W)$.
Then we obviously have
\begin{claim}\label{forms-8}
The pairing $R^k(W)\times R^{n-k}(W)\to M(W)$ is perfect. In
particular one has an isomorphism $R^k(W)\simeq
R^{n-k}(W)^*\otimes M(W)$ commuting with the action of
$Aut_{\HH}W$.
\end{claim}
Remind that the group $Aut_\HH (W)$ of invertible automorphisms of
$W$ is isomorphic to $GL_n(\HH)$, and its complexification is
isomorphic to $GL_{2n}(\CC)$.
\begin{theorem}[Irreducibility property]\label{forms-9}
For $0\leq k\leq n$ the space $R^k(W)\otimes_{\RR}\CC$ is an
irreducible $Aut_{\HH}W$-module. As $GL_{2n}(\CC)$-module it has
highest weight $(\underset{2k \mbox{
times}}{\underbrace{1,\dots,1}},\underset{2(n-k) \mbox{
times}}{\underbrace{0,\dots,0}})$.
\end{theorem}

{\bf Proof.} By Proposition \ref{forms-7} the
$GL_{2n}(\CC)$-module $R^k(W)\otimes_{\RR}\CC$ is multiplicity
free. Let $(\lam_1\geq \dots \geq \lam_{2n})$ be the highest
weight of an irreducible component $H_1$ of this representation
(thus $\lam_i \in \ZZ$). Since it can be realized as a
subrepresentation in the tensor power $Z^{\otimes 2n}$ (where
$Z=Hom_{\chh}(T, \ccw)$) it corresponds to a Young diagram, and
hence $\lam_{i}\geq 0$ for all $i$. By Claim \ref{forms-8} there
exists an irreducible submodule $H_2$ of $R^{n-k}$ such that
$H_1\simeq H_2^*\otimes M(W)$. Let $(\mu_1 \geq \dots \geq
\mu_{2n})$ be the highest weight of $H_2$. Similarly $\mu_i\geq
0$. Then the highest weight of $H_2^*\otimes M(W)$ is equal to
$(1-\mu_{2n},\dots, 1-\mu_1)$. Thus we get
$$\lam_i=1-\mu_{2n-i+1} ,\, i=1,\dots, 2n.$$ Hence $\lam_i\leq 1$
for all $i$. Thus the highest weight of $H_1$ is equal to
$(\underset{l \mbox{
times}}{\underbrace{1,\dots,1}},\underset{2n-l \mbox{
times}}{\underbrace{0,\dots,0}})$.  The center of $GL_{2n}(\CC)$
consists of $\{\lam \cdot Id |\, \lam\in \CC^*\}$, and it acts on
$R^k(W)\otimes_{\RR}\CC$ by multiplication by $\lam ^{2k}$. Hence
$l=2k$. This proves Theorem \ref{forms-9}. \qed

\begin{theorem}\label{forms-10}
The correspondence $W\mapsto R^k(W)$ is a functor from the
category of finite dimensional right $\HH$-modules to the category
of finite dimensional real vector spaces. Namely a morphism
$W_1\to W_2$ induces a morphism $R^k(W_1)\to R^k(W_2)$ in a way
compatible with compositions of morphisms.
\end{theorem}
{\bf Proof.} Let $f:W_1\to W_2$ be a morphism of right
$\HH$-modules. Then $f$ induces a morphism of $\RR$-vector spaces
$f:Sym^k(W_1\otimes \bar W_1)^\sigma\to Sym^k(W_2\otimes \bar
W_2)^\sigma$. We have to show that $f(L^k(W_1))\subset L^k(W_2)$.
Note that $Sym^k(W\otimes \bar W)^\sigma \otimes _{\RR}\CC$ can be
realized (functorially) as a submodule in $Z^{\otimes 2k}$ with
$Z=Hom_{\chh}(T,{}\!^{\textbf{C}} W)$. Then $L^k(W)$ can be
realized as a kernel of corresponding Young symmetrizer which does
not depend on $W$. Hence $f(L^k(W_1))\subset L^k(W_2)$. \qed

Set $R(W):=\oplus_{k=0}^n R^k(W)$ where $n=\dim_{\HH}W$.
\begin{proposition}\label{forms-11}
$R(W)$ is a commutative associative graded algebra with Poincar\'e
duality. Any morphism $f:W_1\to W_2$ induces a homomorphism of
algebras $f:R(W_1)\to R(W_2)$.
\end{proposition}
{\bf Remark.} Recall that the Poincar\'e duality means that
$\dim_{\RR} R^n(W)=1$ and the pairing induced by the
multiplication $R^k(W)\otimes R^{n-k}(W)\to M(W)$ is perfect.

{\bf Proof.} The only property one has to check is
$$L^k(W)\cdot Sym^l(W\otimes_{\HH}\bar W)^\sigma\subset
L^{k+l}(W).$$ Let us check it. Take any $g\in
Sym^{n-k-l}(W\otimes_{\HH}\bar W)^\sigma$. One has $M(L^k(W)\cdot
Sym^l(W\otimes_{\HH}\bar W)^\sigma \cdot g)\subset M(L^k(W)\cdot
Sym^{n-k}(W\otimes_{\HH}\bar W)^\sigma)=0.$ \qed

Now we can define a quaternionic analogue of the space of
translation invariant forms of type $(k,k)$ on a complex space.
Let $V$ be a finite dimensional right $\HH$-module.
\begin{definition}\label{forms-12}
Define $$\Omega^{k,k}(V):= R^k(\bar V^*),$$
$$\Omega^\bullet (V):=\oplus _{k=0}^n\Omega^{k,k}(V).$$
\end{definition}

\begin{theorem}\label{forms-13}
The correspondence $V\mapsto \Omega ^\bullet(V)$ is a
contravariant functor from the category of finite dimensional
right $\HH$-modules to the category of finite dimensional
commutative associative graded algebras. For a fixed $V$ the
graded algebra $\Omega^\bullet (V)$ satisfies the Poincar\'e
duality. For any $k$, $\Omega^{k,k}(V)\otimes_{\RR}\CC$ is an
irreducible $(Aut_{\HH}V\otimes_ \RR \CC)$-module.
\end{theorem}
{\bf Proof.} The proof immediately follows from Theorems
\ref{forms-10}, \ref{forms-9}, and Claim \ref{forms-8}. \qed

Let us denote by $\ch(V)$ the space of hyperhermitian forms on
$V$.
\begin{proposition}\label{forms-15}
There exists canonical isomorphism $\ch(V)=\Omega^{1,1}(V)$.
\end{proposition}
{\bf Proof.} Let $B:V\times V\to \HH$ be a hyperhermitian form. It
defines a map $B_1:V\to Hom_{\HH}(V,\HH)=V^*$. Let us identify $V$
with $\bar V$ as $\RR$-vector spaces. Thus $B_1:\bar V\to V^*$.
Since $B$ is semi-linear with respect to the first argument,
$B_1:\bar V\to V^*$ is a morphism of left $\HH$-modules. Thus
$B_1\in End_{\HH}(\bar V,V^*)=(\bar
V)^*\otimes_{\HH}V^*=\bar{V^*}\otimes _{\HH}V^*.$ The fact that
$B$ satisfies $B(x,y)=\overline{B(y,x)}$ means that $B_1\in
(\bar{V^*}\otimes _{\HH}V^*)^\sigma=\Omega^{1,1}(V)$. Thus the
above construction defines a canonical map $\ch(V)\to
\Omega^{1,1}(V)$. This is an isomorphism. \qed

\begin{proposition}\label{forms-14}
(1) Let $V$ be a right $\HH$-module, $n=\dim_{\HH}V$. There exists
a natural isomorphism
$$(\Omega^{n,n}(V))^{\otimes 2}=\wedge^{4n}_{\RR} (V^*).$$

(2) Let $0\to W\to V\to U\to 0$ be a short exact sequence of
finite dimensional right $\HH$-modules of quaternionic dimensions
$k,n,l$ respectively. Then there exists canonical isomorphism
$$\Omega^{n,n}(V)=\Omega^{k,k}(W)\otimes_{\RR} \Omega ^{l,l}(U).$$
\end{proposition}
{\bf Proof.} (1) Let us fix an $\HH$-basis $e_1,\dots,e_n$ of $V$.
Let $Id\in \ch_n\simeq \ch(V)$ be the identity matrix. Then
$(Id)^n\in \Omega^{n,n}(V)$ spans $\Omega^{n,n}(V)$ over $\RR$.
Let us define isomorphism $\Omega^{n,n}(V)\tilde\to
\wedge_\RR^{4n}V^*$ such that $(Id)^n$ is mapped to
$\wedge_{p=1}^n \left( e_p^*\wedge(e_p\cdot I)^*\wedge(e_p\cdot
J)^*\wedge(e_p\cdot K)^*\right)$ where $*$ denotes taking the
bi-dual basis. It is easy to see that this isomorphism is
independent of a choice of a basis $e_1,\dots,e_n$ of $V$.

(2) Let us define a map $\Omega^{k,k}(W)\otimes_{\RR} \Omega
^{l,l}(U)\to\Omega^{n,n}(V)$ by $x\otimes y\mapsto x\cdot y$. It
is easy to see that this map is well defined and is an
isomorphism. \qed

\begin{remark}\label{forms-3} Let $\dim_\HH V=n$.
The real line $\Omega^{n,n}(V)$ is canonically oriented. The
orientation is defined as follows. Let us fix an arbitrary
$\HH$-basis $e_1,\dots,e_n$ of $V$. Then $\Omega^{1,1}(V)$ is
identified with $\ch(V)\simeq \ch_n$ by Proposition
\ref{forms-15}. Let $Id\in \ch_n$ be the identity matrix. Then
$(Id)^n$ spans $\Omega^{n,n}(V)$. Let us choose the orientation of
$\Omega^{n,n}(V)$ so that $(Id)^n$ is a positive element. It is
easy to check that this orientation is independent of a choice of
a basis of $V$.
\end{remark}

\subsection{Positive forms.}\label{posforms} In this subsection we will define
convex cones of weakly and strongly positive forms in
$\Omega^{k,k}(V)$. The exposition is analogous to the complex case
as in Harvey \cite{harvey} (see also Lelong \cite{lelong}).

Let $V$ be finite dimensional right $\HH$-module, $\dim_{\HH}V=n$.
First recall that by Remark \ref{forms-3} the space
$\Omega^{n,n}(V)$ is a one dimensional real vector space with
canonical orientation. Thus we can define non-negative elements in
$\Omega^{n,n}(V)$. This half line will be denoted by
$\Omega^{n,n}(V)_{\geq 0}$.
\begin{definition}\label{pos-1}
(1) An element $\eta\in \Omega^{k,k}(V)$ is called {\itshape
strongly positive} if it can be presented as a finite sum of
elements of the form $f^*\xi$ where $f:V\to U$ is a morphism of
right $\HH$-modules, $\dim_{\HH}U=k$, $\xi\in
\Omega^{k,k}(U)_{\geq 0}$.

(2) An element $\eta\in \Omega^{k,k}(V)$ is called {\itshape
weakly positive} (or just positive) if for any strongly positive
element $\zeta\in \Omega^{n-k,n-k}(V)$ the product $\eta\cdot
\zeta \in \Omega^{n,n}(V)_{\geq 0}.$
\end{definition}
Clearly positive and weakly positive elements form convex cones.
Let us denote by $C^k(V)$ (resp. $K^k(V)$) the cone of strongly
(resp. weakly) positive elements in $\Omega^{k,k}(V)$.

\begin{proposition}\label{pos-2}
Any strongly positive element of $\Omega^{k,k}(V)$ is weakly
positive, namely $C^k(V)\subset K^k(V)$.
\end{proposition}
{\bf Proof.} It is sufficient to prove the following statement.
Let $V=U\oplus W$, $\dim_{\HH}U=k,\, \dim_{\HH}W=l$. Let $p_U:V\to
U,\, p_W:V\to W$ be the corresponding projections. Let $\eta\in
\Omega^{k,k}(U)_{\geq 0},\, \zeta\in \Omega^{l,l}(V)_{\geq 0}$.
Then one should check that $p_U^*\eta \cdot p_W^*\zeta \in
\Omega^{n,n}(V)_{\geq 0}.$ This is clear. \qed

The proof of the next proposition is obvious.
\begin{proposition}\label{pos-3}
$$C^k(V)\cdot C^l(V)\subset C^{k+j}(V).$$
\end{proposition}

The complex version of the next result is due to Lelong
\cite{lelong}.
\begin{proposition}\label{pos-4}
For $k=0,1,n-1,n$ the cones of weakly and strongly positive
elements coincide.
\end{proposition}
{\bf Proof.} For $k=0,n$ the statement is obvious. Let us prove it
for $k=1$. We may assume that $V=\HH^n$. Let us fix on $\HH^n$ the
standard hyperhermitian form $\sum_{i=1}^n |q_i|^2$. Then as we
have noticed in Subsection \ref{forms} the space
$W\otimes_{\HH}\bar W$ ($W=\bar V^*$) can be identified with the
space of hyperhermitian matrices $\ch_n$. Let us check that
$L^1(W)=0$. It is enough to show that for any $A\in \ch_n, \, A\ne
0$ there exists $B \in \ch_n$ such that the mixed Moore
determinant $det(A,B,\dots, B)\ne 0$. By Claim \ref{det-11} we may
assume that $A$ is diagonal. Say $A=diag(t_1,\dots, t_n)$, $t_1
\ne 0$. Take $B=diag(0,1\dots, 1)$. Then $det(A,B,\dots, B)=t_1\ne
0$. Thus $\Omega^{1,1}(V)$ is identified with $\ch_n$. Let $A\in
K^1(V)$. We may assume that $A$ is diagonal,
$A=diag(t_1,\dots,t_n)$. It is easy to see that $t_i\geq 0$ for
any $i$. This implies that $A\in C^1(V)$. Hence $C^1(V)=K^1(V)$.

Let us assume that $k=n-1$. We may and will assume that $V=\HH^n$.
Let us fix on $\HH^n$ the standard hyperhermitian form
$\sum_{i=1}^n |q_i|^2$. This form gives an identification
$\Omega^{n,n}(\HH^n)\simeq \RR$. By Theorem \ref{forms-13} we have
$\Omega^{n-1,n-1}(\HH^n)=(\Omega^{1,1}(\HH^n))^*\otimes
\Omega^{n,n}(\HH^n)=(\Omega^{1,1}(\HH^n))^*=(\ch_n)^*$. However on
$\ch_n$ we have a non-degenerate bilinear form
$$(X,Y):=Tr(X\cdot Y)$$
where for a quaternionic matrix $Z=(z_{ij})$ by definition
$Tr(Z):=Re \sum_{i} z_{ii}$. Thus using this form one can identify
$\ch_n^*=\ch_n$. One can easily see that with this identification
the cone $\cc\subset \ch_n$ of non-negative definite
hyperhermitian matrices is self-dual, namely $\cc=\{X\in \cc \,|
\, Tr(X\cdot Y)\geq 0 \, \forall Y\in \cc \}$.

Thus we have an identification
\begin{eqnarray}\label{pos-id}
\Omega^{n-1,n-1}(\HH^n)=\ch_n.
\end{eqnarray}
It is easy to see that an element $\eta\in
\Omega^{n-1,n-1}(\HH^n)$ is strongly positive if and only if it
can be represented as a sum with non-negative coefficients of
elements of the form $p^*(vol^{1/2})$ where $p:\HH^n\to L$ is an
{\itshape orthogonal} projection onto a quaternionic subspace $L$
of quaternionic dimension $n-1$, and $vol^{1/2}\in \Omega
^{n-1,n-1}(L)$ is the positive square root of the Lebesgue measure
on $L$ induced by the metric. Let us fix such $L$. Making an
orthogonal transformation we may assume that $L$ is spanned by the
last $n-1$ coordinates. We claim that under the above
identification (\ref{pos-id}) the element $\eta$ corresponds to
the matrix $\frac{1}{n}S$ where $S=\left[ \begin{array}{cccc}
                              1&0&\dots&0\\
                              0&0&\dots&0 \\
                              \multicolumn{4}{c}{\dotfill} \\
                              0&0&\dots&0
                              \end{array}\right]$.
Let us check it. We have to show that for any $X=(x_{ij})\in
\ch_n$ one has $$\det(X,p^*(I_{n-1}),\dots,
p^*(I_{n-1}))=Tr(X\cdot \frac{1}{n}S).$$ Note first that
$vol^{1/2}=(I_{n-1})^{n-1}\in \Omega ^{n-1,n-1}(\HH^n)$ where
$I_{n-1}$ is the identity matrix of size $n-1$. Obviously
\begin{eqnarray}\label{pos-obvi}
p^*(I_{n-1})=\left[
\begin{array}{c|c}
                                                 0&0\\ \hline
                                                 0&I_{n-1}
                                                 \end{array}
                                                 \right].
\end{eqnarray}
We have \begin{eqnarray*} \det(X,p^*(I_{n-1}),\dots,
p^*(I_{n-1}))&=& \frac{1}{n!}\frac{\pt^n}{\pt \lam_1 \dots \pt
\lam_n}\det(\lam_1
X+\sum_{i=2}^n\lam_ip^*I_{n-1})=\\
&&\frac{1}{n!}\frac{\pt^n}{\pt \lam_1 \dots \pt \lam_n}\det(\lam_1
X+(\sum_{i=2}^n\lam_i)p^*I_{n-1}).\end{eqnarray*} Using the
identities (\ref{det-form}) and (\ref{pos-obvi}) we easily get
\begin{eqnarray*}\det(X,p^*(I_{n-1}),\dots,
p^*(I_{n-1}))&=&\frac{1}{n!}\frac{\pt^n}{\pt \lam_1 \dots \pt
\lam_n} \sum_{I\subset \{2,\dots,n\}} (\sum_{i=2}^n\lam_i)^{|I|}
\det M_I(\lam_1X)=\\
\frac{1}{n!}\frac{\pt^n}{\pt \lam_1 \dots \pt
\lam_n}(\sum_{i=2}^n\lam_i)^{n-1}\lam_1x_{11}&=&\frac{1}{n}x_{11}=Tr(X\cdot
\frac{1}{n}S).
\end{eqnarray*}

Hence matrices of the form $p^*(vol^{1/2})$, where $p$ is an
orthogonal projection onto a quaternionic subspace of quaternionic
dimension $n-1$,  are precisely (hyperhermitian) matrices which
can be written at some orthonormal basis is the form $\left[
\begin{array}{cccc}
                              1&0&\dots&0\\
                              0&0&\dots&0 \\
                              \multicolumn{4}{c}{\dotfill} \\
                              0&0&\dots&0
                              \end{array}\right]$.
The set of linear combinations of such matrices with non-negative
coefficients coincides with the cone $\cc$ on non-negative
definite matrices. Hence under the identification (\ref{pos-id})
the cone $C^{n-1}(\HH^n)\subset \Omega^{n-1,n-1}(\HH^n)$
corresponds to the cone $\cc\subset \ch_n$. As we have mentioned
the cone $\cc$ is self-dual. Hence
$K^{n-1}(\HH^n)=C^{n-1}(\HH^n)$. This proves Proposition
\ref{pos-4}. \qed

\begin{proposition}\label{pos-5}
The cones of strongly and weakly positive elements in
$\Omega^{k,k}(V)$ have non-empty interiors.
\end{proposition}
{\bf Proof.} Since $C^k(V)\subset K^k(V)$ by Proposition
\ref{pos-2}, it is enough to check that the cone $C^k(V)$ of
strongly positive elements in $\Omega^{k,k}(V)$ has non-empty
interior. We will use the following well known and simple fact.
\begin{lemma}\label{pos-6}
Let $X$ be a finite dimensional real vector space. Let $C$ be a
convex cone in $X$. Then $C$ has a non-empty interior if and only
if $C-C=X$ where $C-C:=\{c_1-c_2 \, | \, c_1,c_2\in C\}$.
\end{lemma}
Observe first that the cone $C^k(V)$ is $End_\HH(V)$-invariant (as
well as the cone $K^k(V)$). Note that $C^k(V)-C^k(V)$ is a
non-zero $End_{\HH}V$-invariant subspace of $\Omega^{k,k}(V)$. By
the irreducibility property (Theorem \ref{forms-9}) this subspace
must be equal to $\Omega^{k,k}(V)$. Hence by Lemma \ref{pos-6}
$C^k(V)$ has non-empty interior. \qed

\section{More quaternionic pluripotential theory.}\label{psh}
\subsection{Positive currents.}\label{poscur} In this subsection we construct a
quaternionic analogue of form valued measures (currents)
$dd^cf_1\wedge\dots \wedge dd^cf_k$ where $f_i$ are continuous psh
functions.

Let us start with some well known general remarks (see e.g.
\cite{guillemin-sternberg}). Let $M$ be a smooth manifold. Let $E$
be a finite dimensional vector bundle over $M$. Let
$C^\infty_c(M,E)$ denotes the Schwartz space of infinitely smooth
compactly supported sections of $E$ equipped with the standard
topology. Let $|\omega_M|$ denote the line bundle of densities on
$M$.

\begin{definition}\label{cur-1}
A generalized section of $E$ is a continuous linear functional on
$C^\infty_c(M,E^*\otimes |\omega_M|)$.
\end{definition}
The space of generalized sections is denoted by
$C^{-\infty}(M,E)$. Recall also that $C^{\infty}(M,E)\subset
C^{-\infty}(M,E)$. Indeed any $f\in C^{\infty}(M,E)$ defines a
continuous linear functional on $C^\infty_c(M,E^*\otimes
|\omega_M|)$ by $<f,\phi>=\int_M \phi(f)$.

Let now $X$ be a finite dimensional real vector space. The space
of smooth (resp. generalized) functions on $M$ with values in $X$
is denoted by $C^\infty(M,X)$ (resp. $C^{-\infty}(M,X)$).

Let $K\subset X$ be a closed convex cone in $X$. Recall that the
polar of $K$ is denoted by
$$K^\circ :=\{y\in X^* |\, <y,x>\geq 0 \,\, \forall x\in K\}.$$
Then $K^\circ$ is a closed convex cone in $X^*$.
\begin{definition}\label{cur-2}
A {\itshape $K^\circ$-valued smooth density} on $M$ is a section
$\mu\in C^\infty(M,X^*\otimes |\omega_M|)$ such that its value at
any point $p\in M$ has the form
$$\mu(p)=y\otimes l$$
where $y\in K^\circ ,\, l\in |\omega_M|_p,\, l\geq 0$.
\end{definition}
It is easy to see that $K^\circ$-valued smooth densities form a
convex cone in $C^\infty(M,X^*\otimes|\omega_M|)$.
\begin{definition}\label{cur-3}
A generalized $K$-valued function on $M$ is an $\RR$-linear
continuous functional $$f:C^\infty_c(M,X^*\otimes |\omega_M|)\to
\RR$$ such that $f(\phi)\geq 0$ for any smooth compactly supported
$K^\circ$-valued density $\phi$ on $M$.
\end{definition}
\begin{proposition}\label{cur-4}
Assume that $K^\circ \subset X^*$ has non-empty interior. Let $f$
be a generalized $K$-valued function on  $M$. Then $f$ has order
0, i.e.  $f$ is a continuous $\RR$-linear functional on
$C_c(M,X^*\otimes |\omega _M|)$ which takes non-negative values on
$K^\circ$-valued smooth densities.
\end{proposition}
{\bf Proof.} Let us fix a smooth density $m\in C^\infty
(M,|\omega_M|)$ which is strictly positive at each point. Let $f$
be a generalized $K$-valued function. Fix a compact subset $A$ of
$M$. We have to show that there exists a constant $C_A$ such that
for any $\phi\in C^\infty_c(M,X^*\otimes |\omega_M|)$ with $supp
(\phi)\subset A$ one has $$|<f,\phi>|\leq C_A||\phi||_0.$$ Let us
fix a function $\gamma\in C^\infty_c(M,\RR)$ which is equal to 1
on $A$ and $\gamma \geq 0$ on $M$. Let us fix a vector $\xi$ from
the interior on $K^\circ$. There exists a constant $C'_A$ such
that for any $\phi\in C^\infty_c(X,X^*\otimes |\omega_M|)$ with
$supp \phi \subset A$ the function $\psi :=\phi+ C'_A ||\phi||_0
\xi \cdot m\cdot \gamma$ is a $K^\circ$-valued density. Thus
$$f(\phi)=f(\psi - C'_A ||\phi||_0 \xi
\cdot m \cdot \gamma)\geq -||\phi||_0f(C'_A \xi\cdot m\cdot \gamma
)=-C_A||\phi||_0.$$ Replacing $\phi$ by $-\phi$ we obtain the
inverse inequality, and hence
$$|f(\phi)|\leq C_A||\phi||_0.$$ \qed

We will apply the above constructions in the following situation.
Let $V$ be a right $\HH$-module, $\dim_{\HH}V=n$. Fix an integer
$0\leq k\leq n$. Take $X=\Omega^{k,k}(V)$. Then by Claim
\ref{forms-8} one can canonically identify $X^*$ with
$\Omega^{n-k,n-k}(V)\otimes(\Omega^{n,n}(V))^*$. Let us take the
convex cone $K\subset X=\Omega^{k,k}(V)$ to be the cone of weakly
positive elements (in sense of Definition \ref{pos-1}). Recall
that $K$ and $K^\circ$ have non-empty interiors by Proposition
\ref{pos-5}.

Also we will call by {\itshape currents} generalized functions on
$\co$ with values in $\Omega^{\bullet}(V)$.
\begin{definition}\label{cur-5}
Let $\co\subset V$ be an open subset. A current $f\in
C^{-\infty}(\co,\Omega^{k,k}(V))$ is called {\itshape (weakly)
positive} if $f$ is a generalized $K$-valued function on $\co$ in
sense of Definition \ref{cur-1}.
\end{definition}
 It follows from Proposition \ref{cur-4} that any positive current has
order zero. Now we are going to construct a quaternionic analogue
of the expression $dd^c f$. We construct a linear differential
operator $$D_2: C^2(\co, \RR)\to C(\co, \Omega^{1,1})=C(\co,
\ch(V)).$$ Fix any point $v\in V$. The Hessian $B:=Hess_v F$ of
$f$ at $v$ defines a quadratic form on $\rv$. Let us extend it to
$V \otimes_{\RR} \HH$ as a hyperhermitian form $B'$ by the rule
$B'(x\otimes q_1,y\otimes q_2)=\bar q_1B(x,y)q_2$. Define the
subspace
$$V':=\{X-X\cdot I\otimes i-X\cdot J\otimes j-X\cdot K\otimes k\,
| \, X\in V\}\subset V\otimes_\RR \HH.$$ Then $V$ and $V'$ are
naturally isomorphic as right $\HH$-modules.

Consider the restriction of $B'$ to $V'\subset V\otimes_{\RR}
\HH$. Then $B'|_{V'}$ is a hyperhermitian form on $V'$, and in
coordinates it is given by the matrix $\left( \frac{\pt^2 f}{\pt\bar
q_i \pt q_j}\right)$. This defines the map we need
$$D_2: C^2(\co,\RR)\to C(V,\Omega^{1,1}(V)).$$ Note that if $f\in
C^2(\co)\cap P(V)$ then $D_2 f\in \Omega^{1,1}(V)$ is non-negative
pointwise. Hence for such $f$, $D_2 f$ is a  non-negative current.

Now if $k\leq n$ and $f_1,\dots, f_k\in C^2(V,\RR)$ then we can
consider $D_2f_1\cdot \dots \cdot D_2f_k\in C(\co,
\Omega^{k,k}(V))$. Note also that by Proposition \ref{pos-3} this
is a non-negative current provided all $f_i$ are psh. We want to
define this expression as a current for all $f_i\in C(\co)\cap
P(\co)$.
\begin{theorem}\label{cont}
Let $\co\subset V$ be an open subset. For any $k$-tuple of
functions $f_1,\dots, f_k\in C(\co)\cap P(\co)$ one can define a
non-negative $\Omega^{k,k}(V)$-valued current denoted by
$D_2f_1\cdot \dots \cdot D_2f_k$ which is uniquely characterized
by the following properties:

(1) if $f_1,\dots,f_k \in C^2(\co)$ then this current coincides
with the defined above;

(2) if sequences $\{f_i^{(N)}\}\subset C(\co)\cap P(\co)$,
$i=1,\dots,k$, converge to $\{f_i\}$ in $C^0$-topology, i.e.
uniformly on compact subsets of $\co$, then
$$D_2f_1^{(N)}\cdot \dots \cdot D_2 f_k^{(N)}\to D_2f_1\cdot \dots \cdot
D_2f_k \mbox{ as } N\to \infty $$ weakly, i.e. in
$(C^0)^*$-topology.
\end{theorem}
{\bf Proof.} In \cite{alesker-bsm} we have defined a non-negative
measure $(D_2f)^n$ for any function $f\in C(\co)\cap P(\co)$ (in
the notation of \cite{alesker-bsm} $(D_2f)^n$ is identified with
$det\left(\frac{\pt^2 f}{\pt\bar q_i\pt q_j}\right)$) such that
if a sequence $\{f^{(N)}\}\subset C(\co)\cap P(\co)$ converges to
$f$ in $C^0(\co)$-topology then $(D_2f^{(N)})^n\to (D_2f)^n$
weakly. Linearizing the expression $(D_2f)^n$ we define
non-negative measures $D_2f_1\cdot \dots \cdot D_2f_n$ for any
$n$-tuple $f_1,\dots,f_n\in C(\co)\cap P(\co)$, and it satisfies
the same continuity property.

Let $1\leq k<n$. We have to define non-negative currents
$D_2f_1\cdot \dots \cdot D_2f_k$ for $f_1,\dots,f_k\in C(\co)\cap
P(\co)$. Let us fix a basis in $V$. Let us fix $\ch(V)$-valued
functions $A_1,\dots, A_{n-k}\in C(\co,\ch(V))$. Note also that we
can and will identify $\ch(V)$ with the space of hyperhermitian
matrices $\ch_n$. Let us denote for brevity by $\pt^2 u$ the
matrix $\left(\frac{\pt^2 u}{\pt\bar q_i\pt q_j}\right)$. We have
to define measures $det(\pt^2 f_1,\dots, \pt^2 f_k,A_1,\dots,
A_{n-k})$ on $V$ which have the obvious meaning for $f_i \in
C^2(V)$. We will do it by induction in $n$. For $n=1$ the measure
coincides with the usual Laplacian $\Delta f_1$. Assume that
$n>1$. For $f_i\in C^2(\co)$ the expression $det(\pt^2 f_1,\dots,
\pt^2 f_k,A_1,\dots, A_{n-k})$ is linear with respect to each
$A_j$. But every $\ch_n$-valued function can be presented as a
finite linear combination of matrices such that in appropriate
coordinate system each of the summand has the form $\phi \cdot
\left[\begin{array}{cccc}
                      0&\dots&0&0\\
                       &\ddots& & \\
                       0&\dots&0&1
                     \end{array}\right]$ where $\phi\in
                     C(\co,\RR)$.
Thus we may assume that $A_1$ has such a form. Then for any
$f_1,\dots,f_k\in C^2(\co)$ one has by Lemma 1.2.16
\begin{eqnarray*}
\det(\pt^2f_1,\dots,\pt ^2f_k,A_1,\dots,A_k)&=& \\
\frac{1}{n} \phi \cdot
det_{n-1}(M_{\{n\}}(\pt^2f_1),\dots,M_{\{n\}}(\pt^2f_k),
M_{\{n\}}(A_2),\dots, M_{\{n\}}(A_{n-k}))
\end{eqnarray*}
where as previously $M_{\{n\}}(X)$ denotes the minor of the matrix
$X$ obtained by deleting the last row and the last column.
Inductively we are reduced to the following situation. Fix a
decomposition $\HH^n=\HH^k\oplus \HH^{n-k}$. For a $k$-tuple of
functions $f_1,\dots,f_k\in C(\co)\cap P(\co)$ we want to define
measures on $\co$ denoted by $det_k([\pt ^2f_1]_k,\dots,[\pt^2
f_k]_k)$ which depend continuously (with respect to the weak
convergence of measures) in $f_i\in C(\co)\cap P(\co)$. Here
$[X]_k$ denotes the $(k\times k)$-matrix obtained from $X$ by
deleting the last $(n-k)$ rows and columns. For $\psi\in C_c(\co)$
we will define the integral
$$\int_{\HH^n} \psi det_k([\pt ^2f_1]_k,\dots,[\pt^2 f_k]_k)=:
(*)$$
 and show that it is continuous with respect to $f_1,\dots,f_k\in
 C(\co)\cap P(\co)$. Note that if $f_i\in C^2(\co)$ for all $i$
 then
 $$(*)=\int_{\HH^{n-k}}dvol[\int_{\HH^k}\psi det_k([\pt ^2f_1]_k,\dots,[\pt^2
 f_k]_k)dvol].$$
 By Theorem \ref{qpt-2} the inner integral is defined for all $f_1,\dots,f_k\in
 C(\co)\cap P(\co)$ and depends continuously on $f_i$'s. Moreover
 Lemma \ref{det-20} implies the following estimate
 $$|\int_{\HH^k}\psi det_k([\pt ^2f_1]_k,\dots,[\pt^2
 f_k]_k)dvol|\leq\int_{\HH^k}|\psi| det_k([\sum_{i-1}^k\pt
 ^2f_i]_k)dvol.$$
 By Lemma 2.1.9 in \cite{alesker-bsm} for $f\in C(\co)\cap P(\co)$
 one has
 $$\int_{\HH^k}|\psi|det(\pt^2 f) dvol \leq C_k||\psi||_{L^\infty(\HH^n)}
 ||f||^k_{L^\infty(U)}$$
 where $U$ is a compact neighborhood of $supp (\psi)$ and $C_k$ is
 a constant depending on $k$ only. Hence
 $$|\int_{\HH^k}\psi det_k([\pt ^2f_1]_k,\dots,[\pt^2
 f_k]_k)dvol|\leq C'_k ||\psi||_{L^\infty(\HH^n)} (max
 _{i}||f_i||_{C(U)})^k.$$ This estimate and the continuity with
 respect to $f_i$'s of the inner integral in $(*)$  imply that
 $(*)$ is defined for all $f_i\in C(\co)\cap P(\co)$ (since any $f\in C(\co)\cap P(\co)$ can be
 approximated in $C^0$-topology by functions from $C^2(\co)\cap
P(\co)$) and $(*)$ depends continuously on $f_i$'s. \qed
\subsection{Quaternionic Blocki's formula.}\label{qblocki}
Let $\co$ be an open subset in a right $\HH$-vector space $V$,
$\dim_{\HH}V=n$. The complex version of the next result was proved
by Blocki in \cite{blocki}.
\begin{theorem}\label{blocki}
Let $u,v\in P(\co)\cap C(\co)$. Let $2\leq p\leq n$. Then
\begin{equation}\label{blocki-for}
(D_2 \max \{u,v\})^p= (D_2\max\{u,v\})\cdot
\sum_{k=0}^{p-1}(D_2u)^k(D_2v)^{p-1-k} -\sum_{k=1}^{p-1}
(D_2u)^k(D_2v)^{p-k}.\end{equation}
\end{theorem}
{\bf Proof.} First consider the case $p=2$. By continuity of both
sides in (\ref{blocki-for}) we may assume that $u,v$ are smooth.
Let $\chi:\RR \to [0,\infty)$ be a smooth function such that
$\chi(x)=0$ if $x\leq -1$, $\chi(x)=x$ if $x\geq 1$, and $0\leq
\chi ' \leq 1,\, \chi '' \geq 0$ everywhere. Define
$$\psi_j:=v+ \frac{1}{j}\chi(j(u-v)),$$
$$\alpha:=u-v,$$
$$w:=\max\{u,v\}.$$ It is easy to see that $\psi_j\downarrow w$
uniformly on compact subsets and monotonically as $j\to \infty$.
\begin{lemma}\label{308}
$$\left(\frac{\chi(j\alpha)}{j}\right)_{\bar p q}= \chi
'(j\alpha)\cdot \alpha_{\bar p q}+ j\chi ''(j\alpha)\alpha_{\bar
p}\alpha_q.$$
\end{lemma}
{\bf Proof.} $$\left(\frac{\chi(j\alpha)}{j}\right)_{\bar p
q}=\frac{1}{j} \sum_{l,m=0}^3 e_l(\chi (j\alpha))_{x_p^l
x_q^m}\bar e_m=$$
$$\sum_{l,m=0}^3 e_l\left( \chi '(j\alpha)\cdot \alpha_{x_p^l}\right)
_{x_q^m}\bar e_m = \sum_{l,m=0}^3 e_l \left(j \chi ''(j\alpha)
\cdot \alpha _{x_p^l}\alpha _{x_q^m} + \chi '(j\alpha )
\alpha_{x_p^l x_q^m}\right) \bar e_m=$$
$$\chi '(j\alpha) \alpha_{\bar p q} + j \chi ''(j\alpha)\alpha
_{\bar p}\alpha _q.$$ \qed

Thus from Lemma \ref{308} we obtain
$$ (\psi_j)_{\bar p q}=v_{\bar p q}+ \chi' (j\alpha)(u-v)_{\bar p q}+
j \chi ''(j\alpha)\alpha _{\bar p}\alpha _q=$$
$$\chi '(j\alpha)u_{\bar p q}+ (1-\chi'(j\alpha ))v_{\bar p q}+j \chi ''(j\alpha)\alpha
_{\bar p}\alpha _q.$$ Since $0\leq \chi'\leq 1$ and $\chi''\geq 0$
this implies that $\psi_j$ is psh. From the definition of $\psi_j$
we have
\begin{equation}\label{2}(D_2
\psi_j)^2=(D_2 v)^2+ 2(D_2
v)\left(D_2\left(\frac{\chi(j\alpha)}{j}\right)
\right)+\left(D_2\left( \frac{\chi(j\alpha)}{j}\right) \right)^2.
\end{equation}
We have weak convergence
\begin{eqnarray}
& & (D_2\psi_j)^2\to (D_2 w)^2,\label{3}\\
& & D_2v \cdot D_2\left(\frac{\chi(j\alpha)}{j}\right)\to
D_2(w-u)\cdot D_2 v \label{4}.
\end{eqnarray}
Let us study the last term in (\ref{2}), namely $\left(D_2\left(
\frac{\chi(j\alpha)}{j}\right) \right)^2$.
\begin{lemma}\label{309}
Let $A= \left[ \begin{array}{c}
               a_1\\
               \vdots\\
               a_n
               \end{array}\right]\in \HH^n.$
Then in $\Omega^{2,2}(\HH^n)$ one has $(AA^*)^2=0$.
\end{lemma}

{\bf Proof.} One has to check that for all $X_1,\dots, X_{n-2}\in
\ch_n$
$$\det(AA^*,AA^*,X_1,\dots ,X_{n-2})=0.$$ Multiplying $A$ by an
invertible matrix one may assume that $A= \left[ \begin{array}{c}
               1\\
               0\\
               \vdots\\
               0
               \end{array}\right].$
Then $AA^*=\left[ \begin{array}{cccc}
                     1&0&\dots&0\\
                     0&0&\dots&0\\
                     \multicolumn{4}{c}\dotfill\\
                     0&0&\dots&0\\
                     \end{array} \right].$
The result now follows from Lemma \ref{det-19}. \qed

From Lemmas \ref{308} and \ref{309} one gets in $\Omega^{2,2}(V)$:
\begin{equation}\label{star}
\left(\left( \frac{\chi(j\alpha)}{j}\right)_{\bar p q} \right)^2=
\left( \chi '(j\alpha)^2 (\alpha_{\bar p q})+ 2j\chi
'(j\alpha)\chi''(j\alpha)(\alpha_{\bar p}\alpha_q)\right) \cdot
(\alpha_{\bar p q}).
\end{equation}

Let $\gamma:\RR \to \RR$ be such that $\gamma' =(\chi')^2$. Then
we have
\begin{lemma}\label{3010}
$$\left(D_2\left( \frac{\chi(j\alpha)}{j}\right) \right)^2=
D_2\left(\frac{\gamma(j\alpha)}{j}\right)\cdot D_2\alpha.$$
\end{lemma}
Let us postpone the proof of Lemma \ref{3010}, and finish the
proof of Theorem \ref{blocki} for $p=2$. One can choose $\gamma$
so that $\gamma (-1)=0$. Then $\frac{\gamma(jx)}{j}\downarrow
\max\{0,x\}$ uniformly on compact subsets and monotonically as
$j\to \infty$, and
$$\left(D_2\left( \frac{\chi(j\alpha)}{j}\right) \right)^2\to
D_2(w-v)\cdot D_2\alpha \mbox{ weakly }.$$ This and (\ref{2}),
(\ref{3}), (\ref{4}) imply
\begin{eqnarray*}(D_2w)^2&=&(D_2v)^2
+2D_2(w-v)\cdot D_2v+D_2(w-v)\cdot D_2(u-v)\\
&=& D_2w\cdot (D_2u+D_2v)-D_2u\cdot D_2v.
\end{eqnarray*}
This implies Theorem \ref{blocki} for $p=2$. It remains to prove
Lemma \ref{3010}.

{\bf Proof of Lemma \ref{3010}.} We have
\begin{eqnarray*}
\left(\frac{\gamma(j\alpha)}{j}\right)_{\bar p q}&=&
\gamma'(j\alpha)\alpha_{\bar p q}+\alp_{\bar p} \cdot
(\gamma'(j\alp))_q\\
&= & (\chi'(j\alp))^2\alpha_{\bar p q}+ 2j \chi'(j\alp)
\chi''(j\alp)\alp_{\bar p}\cdot \alp_q.
\end{eqnarray*}
This and (\ref{star}) imply Lemma \ref{3010}. \qed

It remains to prove Theorem \ref{blocki} for $p>2$. Set
\begin{equation}\label{xyz}
x:=D_2u,\, y:=D_2u,\, z:=D_2\max\{u,v\}.
\end{equation}
One has
\begin{equation}\label{form}
z^2=z(x+y)-xy.
\end{equation}

One has to show that
$$z^p=z\sum_{k=0}^{p-1} x^ky^{p-1-k}-\sum_{k=1}^{p-1}x^ky^{p-k}.$$
Let us prove it by induction in $p$. For $p=2$ it is already
proved. Assume that the statement is proved for $p$. Let us prove
it for $p+1$. By the assumption of induction one has
\begin{eqnarray*}
z^{p+1}&=& z(z\sum_{k=0}^{p-1}
x^ky^{p-1-k}-\sum_{k=1}^{p-1}x^ky^{p-k})\\
&=& (z(x+y)-x y)\sum_{k=0}^{p-1} x^ky^{p-1-k}-
z\sum_{k=1}^{p-1}x^ky^{p-k}\\
&=& z\sum_{k=0}^{p} x^ky^{p-k}-\sum_{k=1}^p x^ky^{p+1-k}.
\end{eqnarray*}
\qed
\begin{theorem}\label{val}
Let $u,v\in C(\Omega)\cap P(\Ome)$. Assume that $\min \{u,v\}\in
P(\Ome)$. Then for $1\leq p\leq n$ one has
$$(D_2u)^p+(D_2v)^p= (D_2 \min\{u,v\})^p+(D_2 \max\{u,v\})^p.$$
\end{theorem}
{\bf Proof.} Note that $\min\{u,v\}=u+v-\max\{u,v\}$. Hence
$$D_2\min\{u,v\}=x+y-z$$
in the notation (\ref{xyz}). Thus we have to show that
$$x^p+y^p=z^p+(x+y-z)^p.$$
Let us prove it be the induction in $p$. For $p=1$ the statement
is clear. Let us assume that the statement is true for $p-1$. Let
us prove it for $p$. We have
\begin{eqnarray*}
(x+y-z)^p&=& (x+y-z)(x^{p-1}+y^{p-1}-z^{p-1})\\
&=& (x^p+y^p-z^p)+[2z^p+
xy^{p-1}-xz^{p-1}+yx^{p-1}-yz^{p-1}-zx^{p-1}-zy^{p-1}].
\end{eqnarray*}
Let us denote the summand in square brackets by $A$. By Theorem
\ref{blocki} and the assumption of induction we have
\begin{eqnarray*}
A&=& 2[z\sum_{k=0}^{p-1}
x^ky^{p-1-k}-\sum_{k=1}^{p-1}x^ky^{p-k}]+xy^{p-1}-\\
& &- (x+y)[z\sum_{k=0}^{p-2}x^ky^{p-2-k}
-\sum_{k=1}^{p-2}x^ky^{p-1-k}]
+yx^{p-1}-z(x^{p-1}+y^{p-1})\\
&=&z\{2\sum_{k=0}^{p-1}x^ky^{p-1-k}-
(x+y)\sum_{k=0}^{p-2}x^ky^{p-2-k}- (x^{p-1}+y^{p-1})\}-\\
& &-2\sum_{k=1}^{p-1}x^ky^{p-k}
+xy^{p-1}+(x+y)\sum_{k=1}^{p-2}x^ky^{p-1-k}+yx^{p-1}=0.
\end{eqnarray*}
This implies Theorem \ref{val}. \qed

\section{Valuations and complex and quaternionic analysis.}\label{valuations}
In this section we discuss some applications of results from
complex and quaternionic pluripotential theory discussed in
Section \ref{psh} to construction of continuous valuations. In
Subsection \ref{kazar} we discuss complex case, namely
Kazarnovskii's pseudovolume and its generalizations. In Subsection
\ref{vq} we discuss two quaternionic versions of Kazarnovskii's
pseudovolume and generalizations.
\subsection{Kazarnovskii's pseudovolume.}\label{kazar}
Let $V=\CC^n$ be a hermitian space. Let us consider on the
Grassmannian $Gr_n(V)$ of {\itshape real} $n$-dimensional
subspaces the  function $f$ such that for any $L\in Gr_n(V)$
$f(L)$ is equal to the coefficient of the area distortion under
the orthogonal projection from $L$ to $iL^\perp$. Sometimes $f(L)$
is denoted as $|\cos(L,iL^\perp)|$ (the absolute value of the
cosine of the angle between $L$ and $iL^\perp$). Note also that
$f$ is proportional to $L\mapsto vol_{2n}(Q_L+ i\cdot Q_L)$ where
$Q_L$ denotes the unit cube in $L$ as previously.

Let $\cp(V)$ denote the class of convex polytopes in $V$. For a
polytope $P\in \cp(V)$ let $\cf_k(P)$ denote the set of
$k$-dimensional faces of $P$. For a face $F\in \cf_k(P)$ let
$\gamma(F)$ denote the measure of the exterior angle at $F$. Let
$\bar F$ denote the only $k$-dimensional linear subspace parallel
to $F$.

Consider the valuation $\phi_f$ on $\cp(\CC^n)$ defined by
$$\phi_f (P)=\sum_{F\in \cf_k(P)} f(\bar F) vol_n(F) \gamma (F).$$
This valuation will be called Kazarnovskii's pseudovolume since it
was introduced and studied by B. Kazarnovskii
\cite{kazarnovskii-81}, \cite{kazarnovskii-84} in connection with
counting of zeros of exponential sums. From the integral geometric
point of view Kazarnovskii's pseudovolume was studied in
\cite{alesker-jdg}. The first part of following result is due to
Kazarnovskii \cite{kazarnovskii-81}, the second part was proved in
\cite{alesker-jdg}.
\begin{theorem}[\cite{kazarnovskii-81}]
The Kazarnovskii pseudovolume admits an extension by continuity to
$\ck(\CC^n)$. This extension is a continuous translation invariant
$n$-homogeneous valuation on $\ck(\CC^n)$ which is
$U(n)$-invariant.
\end{theorem}
We will discuss below the reason why it is true. Kazarnovskii has
given several formulas for the pseudovolume. Let us start with the
following one. For a convex set $\ck(\CC^n)$ let us denote by
$h_K$ its supporting functional. Recall its definition:
$$h_K(x)=\sup_{y\in K}(x,y).$$
\begin{theorem}[\cite{kazarnovskii-81}]
Kazarnovskii's pseudovolume $P(K)$ is equal to
$$\frac{1}{\kappa_n}\int_{B}(dd^c h_K)^n=\frac{2^{2n}n!}{\kappa_n} \int_{B}\det (\frac{\pt^2 h_K}{\pt\bar z_i
\pt  z_j})d vol$$ where $B$ is the unit ball in $\CC^n$ and
$\kappa_n$ is the volume of the unit $n$-dimensional ball.
\end{theorem}

In this theorem the expression under the integral is understood as
a measure (i.e. in the generalized sense). Thus the continuity of
Kazarnovskii's pseudovolume follows from Theorem \ref{ca-3}. We
have also the following result.
\begin{theorem}\label{kazgen}
Let $0\leq k\leq n$.  Fix $\psi\in C_0(\CC^n, \Omega^{n-k,n-k})$.
Then $K\mapsto \int_{\CC^n}(dd^ch_K)^k\wedge \psi$ defines a
continuous translation invariant $k$-homogeneous valuation on $\ck
(\CC^n)$.
\end{theorem}

In this paper we will prove a quaternionic analogue of this
result, Theorem \ref{quagen}. The proof of Theorem \ref{kazgen}
can be obtained similarly using the corresponding results from
complex analysis. We omit the details of the proof of this
theorem, and prove instead the quaternionic version.


\subsection{Valuations on quaternionic spaces.}\label{vq} Let $V=\HH^n$ be
the right quaternionic hyperhermitian space. We also fix the
standard hyperhermitian product $(x,y)=\sum_{i=1}^n\bar x_iy_i$.
The next result is a quaternionic version of Theorem \ref{kazgen}.
\begin{theorem}\label{quagen}
Let $0\leq k\leq n$. Fix $\psi\in C_0(V,\Omega^{n-k,n-k}(V)\otimes
\Omega^{n,n}(V))$. Then
$$K\mapsto \int_V (D_2 h_K)^k\cdot \psi$$
is a translation invariant continuous $k$-homogeneous valuation on
$\ck(V)$.
\end{theorem}
{\bf Proof.} Translation invariance is obvious. Continuity follows
from Theorem \ref{cont}. To prove the valuation property let us
observe first that if $K=K_1\cup K_2$ with $K_1,K_2,K\in \ck(V)$
then
$$h_K=\max\{h_{K_1},h_{K_2}\},\, h_{K_1\cap
K_2}=\min\{h_{K_1},h_{K_2}\}.$$ Hence the result follows from
Theorem \ref{val}. \qed

From Theorem \ref{quagen} we immediately get the following
corollary which provides new examples of $Sp(n)Sp(1)$-invariant
valuations on $\HH^n$.
\begin{corollary}\label{invariant}
Let $0\leq k\leq n$. Fix $\psi\in C_0(V,\Omega^{n-k,n-k}(V)\otimes
\Omega^{n,n}(V))$ which is invariant under the group $Sp(n)Sp(1)$.
Then
$$K\mapsto \int_V (D_2 h_K)^k\cdot \psi$$
is a translation invariant $Sp(n)Sp(1)$-invariant continuous
$k$-homogeneous valuation on $\ck(V)$.
\end{corollary}

Let us define a {\itshape quaternionic pseudovolume} $Q$ which
will be a quaternionic version of Kazarnovskii's pseudovolume. Let
us consider on the Grassmannian of {\itshape real} $n$-dimensional
subspaces the following function:
$$f(L)=\sqrt{vol_{4n} (Q_L+Q_L\cdot i+Q_L\cdot j +Q_L \cdot k)}$$
where $Q_L$ is the unit cube in $L$ as previously. Consider the
valuation $Q$ on $\cp(\HH^n)$ defined by
$$Q (P)=\sum_{F\in \cf_n(P)} f(\bar F) vol_n(F) \gamma (F).$$
\begin{theorem}\label{vq-2}
The valuation $Q$ extends by continuity to $\ck(\HH^n)$. This
extension is a continuous translation invariant $n$-homogeneous
$Sp(n)Sp(1)$-invariant valuation on $\ck(\HH^n)$.  Moreover for
any $P\in \cp(\HH^n)$
\begin{equation}\label{qps}
Q(P)=\frac{1}{\kappa_{3n}} \int_{D}\det (\frac{\pt^2 h_P}{\pt\bar q_i
\pt  q_j}) dvol
\end{equation}
where $D$ is the unit ball in $\HH^n$ and $\kappa_{3n}$ is the
volume of the unit $3n$-dimensional ball.
\end{theorem}
 Clearly it is enough to prove only equality (\ref{qps}). The
rest follows from Theorem \ref{quagen}. In order to prove equality
(\ref{qps}) we will prove a more precise statement. We will
describe explicitly the measure $\det (\frac{\pt^2 h_P}{\pt\bar q_i
\pt q_j}) dvol$.  We need to introduce more notation.

Let us fix a polytope $P\in \cp(V)$. For a face $F\subset P$ let
us define
$$F^\vee:=\{u\in V^*|\, u(v)=h_P(u) \, \forall v\in F\}$$
where $h_P$ is the supporting function of $P$ as previously. One
can easily check the following properties:

(1) $F^\vee$ is a closed convex polyhedral cone;

(2) $\dim F +\dim F^\vee =\dim_\RR V$;

(3) $F_1\subset F_2$ if and only if $F^\vee_1\supset F_2^\vee$;

(4) $int F^\vee \cap intG^\vee = \emptyset$ for $F\ne G$ where
$int$ denotes the relative interior of a cone;

(5) if $ int F^\vee \cap G^\vee\ne \emptyset $ then $F^\vee
\subset G^\vee$;

(6) the union of $int F^\vee$ when $F$ runs over all non empty
faces of $P$ is equal to $V^*$.

\begin{theorem}\label{vq-3}
Let $P\in \cp(\HH^n)$. The measure $\det (\frac{\pt^2 h_P}{\pt\bar q_i
\pt  q_j}) dvol$ is supported on the union of $F^\vee$ where
$F$ runs over all $n$-dimensional faces of $P$. For an
$n$-dimensional face $F$ of $P$ the restriction of $\det
(\frac{\pt^2 h_P}{\pt\bar q_i \pt q_j}) dvol$ to $F^\vee$ is
equal to $\sqrt{vol_{4n} (Q_{F^\vee}+Q_{F^\vee}\cdot
i+Q_{F^\vee}\cdot j +Q_{F^\vee} \cdot k)} \cdot  vol_{F^\vee}$
where $Q_{F^\vee}$ denotes the unit cube in the span of $F^\vee$,
and $vol_{F^\vee}$ denotes the volume form in the span of $F^\vee$
induced by the Euclidean metric on $\HH^n$.
\end{theorem}
{\bf Proof.} Let us denote for brevity by $\mu_P$ the measure
$\det (\frac{\pt^2 h_P}{\pt\bar q_i \pt  q_j}) dvol$. Since the
supporting functional $h_P$ is homogeneous of degree 1 then the
measure $\mu_P$ is homogeneous of degree $3n$, namely for any
compact set $A$ and any $\lam>0$ one has $\mu_P(\lam A)=\lam
^{3n}\mu_P(A)$. Next we will need the following claim.
\begin{claim}\label{vq-4}
Let $F$ be a face of $P$. Assume that $F$ contains 0. There exists
an open subset $U\subset V^*$ containing $int F^\vee$  such that
$h_P|_U$ is invariant under translations with respect to span
$F^\vee$.
\end{claim}
To prove this claim it is enough to observe that one can choose
$U$ to be the interior of $\cup_G G^\vee$ where $G$ runs over all
vertices of $F$.

Let us now fix a face $F$ of $P$. By Claim \ref{vq-4} the
restriction of $\mu_P$ to $int F^\vee$ is a translation invariant
measure on $int F^\vee$. Hence it is proportional to the Lebesgue
measure. We have previously mentioned that $\mu_P$ is homogeneous
of degree $3n$. Hence if $\dim F \ne n$ the restriction of $\mu_P$
to $int F^\vee$ vanishes. Let us now assume that $\dim F =n$. Thus
the restriction of $\mu_P$ to $F^\vee$ is proportional to the
Lebesgue measure $vol_{3n}$ induced by the Euclidean metric.

Let us denote by $U$ the interior of the union $\cup_G G^\vee$
where $G$ runs over all vertices of $F$. Then $U$ is a non empty
open subset of $V^*$, and for any $y\in U$ one has $$h_P(y)=\sup
\{y(x)\, |\, x\in F\} = h_F(y).$$ Hence we may replace $P$ by $F$,
namely we will assume that $\dim P=n$, $F=P$. Let us prove that if
$\dim P=n$ then the measure $\mu_P$ restricted to $P^\perp$ is
equal to $\sqrt{vol(Q_{P^\perp}+Q_{P^\perp}\cdot I
+Q_{P^\perp}\cdot J+Q_{P^\perp}\cdot K)} vol_{P^\perp}$. Let us
fix an orthonormal basis $\{\xi_1,\dots,\xi_n\}$ in $span P$.
\def\tp{\tilde P}
We can choose a polytope $\tp\subset \RR^n\subset \HH^n$ and an
$\HH$-linear operator $$A:\HH^n\to \HH^n$$ such that $A(\tp)=P$
and for the standard basis $e_1,\dots, e_n$ of $\HH^n$
$$A(e_i)=\xi_i,\, i=1,\dots, n.$$
\def\mhp{\frac{\pt^2 h_P}{\pt\bar q_i\pt  q_j}}
\def\mhtp{\frac{\pt^2 h_{\tp}}{\pt\bar q_i\pt q_j}}
\def\ra{{}\!^ {\mathbb{R}} A}
Then clearly $h_P = h_{\tp} \circ A^*$. Hence by Proposition
\ref{qpt-1.5} for any $q\in \HH^{n*}$ one has
$$\det (\mhp (q))=\det (AA^*)\det(\mhtp (A^* q)).$$
We may assume that $A$ is invertible. Let $w=A^* q$. Then we
obtain
$$dw= |\det (\ra)|dq= (\det AA^*)^2 dq$$
where the last equality follows from Proposition \ref{det-18}.
Hence we conclude
$$\det (\mhp (q)) dq= \det(AA^*)^{-1}\det(\mhtp (w)) dw.$$
\begin{lemma}\label{vq-5}
The measure $\det (\mhtp (w))dw $ is equal to the Lebesgue measure
on $\RR^{n\perp}$ induced by the standard Euclidean metric on
$\HH^n$ times $vol(\tilde P)$.
\end{lemma}
Lemma \ref{vq-5} easily follows from the observation that
$\frac{\pt^2 h_{\tp}}{\pt\bar q_i\pt q_j}=\frac{\pt^2
h_{\tp}}{\pt x_i\pt x_j}$ where $(x_1,\dots, x_n)$ are standard
coordinates on $\RR^n$, and the fact that if $h$ is the
restriction of $h_{\tp}$ to $\RR^{n*}$ then $$\det(\frac{\pt^2
h}{\pt x_i\pt x_j}(w))dw=vol(P)\delta _0$$ where $\delta_0$ is the
delta-measure at the origin 0.

Thus we have to compute $(A^*)^{-1}_*(vol_{\RR^{n\perp}})$ as a
measure on $(span P)^\perp$. For any two Euclidean spaces
$L^{(1)},\, L^{(2)}$ of the same dimension and for any linear map
$ \phi:V^{(1)}\to V^{(2)}$ let us denote by $|\det \phi|$ the
coefficient of the volume distortion, namely for any compact set
$X$ of positive measure $|\det \phi|=\frac{vol \phi(X)}{vol X}$.
One can easily prove the following linear algebraic lemma.
\begin{lemma}\label{vq-6}
Let $V^{(1)},\, V^{(2)}$ be Euclidean vector spaces of the same
dimension. Let $$\phi:V^{(1)}\to V^{(2)}$$ be a linear operator.
Let $L^{(1)}\subset V^{(1)}$ be a linear subspace. Denote
$L^{(2)}=\phi(L^{(1)})$. Assume that $\phi |_{L^{(1)}}:L^{(1)}\to
L^{(2)}$ is a linear isomorphism. Set
$$\psi:=\phi^*|_{L^{(2)\perp}}:L^{(2)\perp}\to L^{(1)\perp}.$$
Then $|\det \psi|=\frac{|\det \phi|}{|\det \phi|_{L^{(1)}}|}.$
\end{lemma}
In our situation $L^{(1)}=\RR^n$, $\phi|_{L^{(1)}} $ is an
isometry. Hence by Lemma \ref{vq-6} one has
$$(A^*)^{-1}_*(vol_{\RR^{n\perp}})=|\det (\ra)|\cdot
vol_{P^\perp}=(\det (AA^*))^2 vol_{P^\perp}.$$ Hence
$$\det(\mhp)dq =\det (AA^*) vol_{P^\perp}=\sqrt{\det
(\ra)}vol_{P^\perp}.$$ But $\det (\ra) =vol_{4n} (Q_{span
P}+Q_{span P}\cdot I+ Q_{span P}\cdot J+ Q_{span P}\cdot K)$. This
proves Theorem \ref{vq-3}. \qed

\end{document}